# Helen of Troy, and the birth of Fuzzy Logic

## Kyriakos Papadopoulos [a], P. Tseliou [b], B.K. Papadopoulos [c]

[a] Department of Mathematics, Kuwait University, PO Box 5969, Safat 13060, Kuwait, emai: kyriakos.p@aum.edu.kw

[b] Department of Sciences of Education in Preschool Age, Democritus University of Thrace, Nea Hili 68 100, Alexandroupolis, email: panatsel1@psed.duth.gr

[c] Department of Civil Engineering, Democritus University of Thrace, Vas. Sofias 12, 67 100 Xanthi, email: papadob@civil.duth.gr

ABSTRACT

The poem Helen of the Nobel laureate George Seferis was inspired by the anti-war play Helen of Euripides. In his poem, Seferis empathizes with the hero of the tragedy, Teucer, who opposed the involvement of The Gods in the lives of the humans, posing unanswered and contradictory questions. With the verse "What is god? What is not a god? And what is there in between?" the ancient poet Euripides sets foundations to the kind of logic that one can consider as a predecessor of Fuzzy Logic. It is worth noting that when Seferis received the Nobel Prize in Stockholm in 1963, he stated: "Right now I feel I am a contradiction myself", a sentence charged with the new language of mathematical logic of the 20$^{th}$ century. Interplaying with the famous words of Karl Weierstrass ``…it is true that a mathematician who is not somewhat of a poet, will never be a perfect mathematician'' in this article we discuss, through two poems, of how poetry and mathematical logic might have influenced each other.

## 1. What is Fuzzy Logic About.

Fuzzy Logic is aimed to be the sort of mathematical logic in which all propositions have a degree of truth between 0 and 1, including 0 and 1. For example, by using fuzzy logic one can support that ``a person x is good by a degree of truth 0.6 and bad by a degree of truth of 0.3". It looks like that such an answer fits better in our complicated reality. If, for example, a group of people were asked the question: ``Did you like that movie?", some would respond ``I liked it", some others ``I did not like it" and many others will say ``I liked it so and so". That last answer can be considered as the "most fuzzy" one, since it supports that the movie was good by a degree of 0.5 and, hence, it leads us to dwell upon the region of fuzzy logic, a field which tries to cover the big gap between 0 and 1.

The foundations of Fuzzy Logic were set in the 1965s, by the university of Berkeley professor, Lotfi Zafeh [1].

## 2. Helen of Troy.

In its official announcement in Stockholm, in the year 1963, the Swedish Academy underlined that George Seferis received the Nobel Prize due to the unique lyricism of his writing which is inspired by his deep emotions for the Hellenic civilisation. Seferis addressed a speech [5] for his

---

nomination, where among other things he said: "I feel at this moment that I am a living contradiction." We emphasize the word contradiction, and we will return back to it soon.

In this article, we are concerned with the mathematical structures that one can extract from a particular poem of Seferis, leaving the artistic part of poetry to the experts and lovers of literature. Through "Helen", we will embark on a quest along with Teucer and the poet himself, to find our own truth closely linked to the realm of Fuzzy Logic.

We first start from the Helen of Euripides, and from the particular fragment:

«ὅτι θεὸς ἢ μὴ θεὸς ἢ τὸ μέσον τίς φησ᾽ ἐρευνάσας βροτῶ».

``What mortal claims, by searching to the utmost limit, to have found out the nature of God, or of his opposite, or of that which comes between, seeing as he doth this world of man tossed to and fro by waves of contradiction and strange vicissitudes?'' [2]

and we continue with the Helen of Seferis:

«τ᾽ είναι θεός; τι μη θεός; και τι τ᾽ ανάμεσό τους; »

``What is a god? What is not a god? And what is there in between them?'' [3]

Teucer, the main hero in Euripides' Helen, is an excellent archer who, after the Trojan War and many adventures, ends up in Cyprus. Seferis, on the other hand, used to be a member of the diplomatic mission in Cyprus, during the tough period around 1953. Both our protagonists, Teucer and Seferis, have identical views and philosophize about the truths in life. Teucer fought in the ten year long Trojan War, for a woman, the "pretty Helen". Seferis lived in a tough period when Greece had to fight for freedom during World War II and experienced the Greek civil war and the struggle of the Cypriots for their independence. Both characters conclude that the numerous possibilities of bloodshed or death take back seats as the brains of mankind cloud over, allowing them to see reality (thus, logic too) from a different perspective.

## 3. The denial of Classical Logic from Seferis.

Through the following apparently simple questions but not so, actually "What is a god? What is not a god? And what is there in between them?" Seferis wonders whether the classical logic of Aristotle (the classical logic of true or false, 0 or 1) is valid and introduces mathematical concepts and questions that were only answered years later: If "God" represents 1 and "Not a God" represents 0, is there something in between 0 and 1? Therefore, the poet ponders on finding the answer to the simple question of "why". The traditional, classical mathematical logic, that one of 0 and 1, is effective when used in mathematical expressions but, everyday life has proven to be much more complicated. For example, one could claim that proposition p: "2 is and odd number" is false or that it has a degree of truth equal to 0. Another might suggest that proposition p: "number 4 is a perfect square" is true with a degree of 1. The controversy starts when one attempts to apply classical logic to everyday life; if we allege that 2,000 euros is a high monthly income for an employee, then considering a salary of 1,999 euros to be a not-high one, does not make much sense. This contradiction occurs due to the so called "Law of Contradiction" of classical logic: There cannot be a proposition which is true and false, simultaneously. According to this, we consider classical logic to be a contradictory one, when we try to interpret real life phenomena.

The success of the accuracy of the science of mathematics is based on Aristotle and his predecessors, who tried to formulate a short theory of mathematical logic and mathematical science: "The Laws of Thought". One of the laws, the so called "Law of the Excluded Middle" states that every proposition must be either true or false. The latter was initially proposed by Parmenides (5th century B.C.) instigating the immediate and strong reaction of his "opponents": Heracletus claimed that a proposition can be simultaneously true and false. Plato also disputed the above law, supporting that a third option exists between true and false, when both "fall to exist". However, there are modern mathematicians and philosophers who support the notion that credit for the foundations of fuzzy logic should be attributed to Buddha who preached that good and bad coexist.

In 1920, the Polish mathematician Jean Lucasiewicz described a system of three-valued logic. The third value that he suggested can be described as feasible and is connected with an arithmetic value between "true or false". Later on, he introduced a system of four-valued and a five-valued logic system, stating that nothing could hinder the way to creating a multi-valued logic. He ended up using the four-valued logic system because he found that it was more pertinent to Aristotle's logic. In 1965, Lotfi Zadeh published an innovative paper under the title "Fuzzy Sets" [4] which dwelled upon the mathematical logic of the theory of fuzzy sets and, inevitably, fuzzy logic. This particular logic, is a generalization of Aristotle's logic, the multi-value logic, introducing the membership function as a quota of truth, between the endpoints of the interval [0,1], where 1 corresponds to "true" and *0* to "false".

In some philosophy books, we turn to discussions concerning the above depicted as a line where the face of Aristotle is placed on the two end points [0,1] and Buddha is placed in the middle (1/2). According to this recursively defined collection of strings, Aristotle's logic is the absolute logic of 0 or 1; the logic of true or false, the Seferis' "What is a God? What is not a God". Buddha might not reflect the phrase "And what is there in between them" but it symbolizes the maximum uncertainty, the everyday phrase "so and so". In a way, the picture shown below, with the portrait or Aristotle place over 0 and over 1 and of Buddha over 1/2, mirrors the three-valued logic of Jean Lucasiewicz:

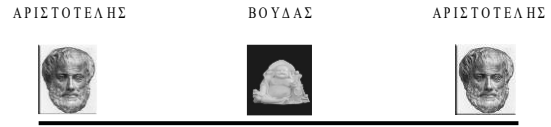

ΑΡΙΣΤΟΤΕΛΗΣ    ΒΟΥΔΑΣ    ΑΡΙΣΤΟΤΕΛΗΣ

**4. Seferis wonders if an extension of classical logic free of contradictions exists.**

At this point, we should note that everything related to classical logic, is interwoven with 0 or 1. The classical set can be defined using 0 or 1 as follows:

We consider a reference set *X*. A subset *A* of *X* ($A \subset X$) is a set with elements of *X*. Every subset *A* of *X* can be described with:

- **Indication:** Indicate all elements of *A*, one by one.
- **Description:** Describe subset *A*, based on element properties.

**Example:** If $X = N$ is the set of natural numbers, the set of even numbers greater than 5 and lower than 7, can be described using the above methods:

**1.** $A = \{6\}$
**2.** $A = \{x \in N | 5 < x < 7\}$

Another way of describing subset *A* of a reference set *X* is by using its characteristic function $\chi_A$. For every subset $A \subset X$, the characteristic function is:

$$x \in X \to \chi_A \in \{0,1\} \; \text{ή} \; \chi_A(x) = \begin{cases} 1, & \text{αν } x \in A, \\ 0, & \text{αν } x \notin A. \end{cases}$$

We can equate subset $A \subset X$ with the characteristic function $\chi_A$ (in algebraic terms: there is an isomorphism from the set of subsets of *X* to the set of characteristic functions).

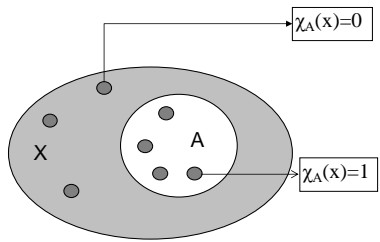

The following example clarifies all of the above. Using classical sets, the set of high temperature can be presented as follows:

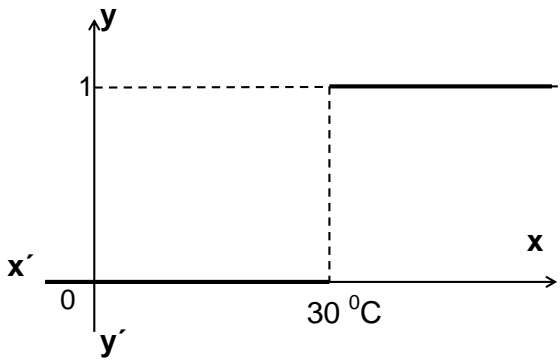

Namely, the classical set of high temperature is given by the function:

$$A(x) = \begin{cases} 1, & \alpha\nu \; x \geq 30 \\ 0, & \alpha\nu \; x < 30 \end{cases}$$

The classical set is interlaced with "1 or 0" or Seferis' "God or Not a God" and is governed by the following laws:

- **Law of Contradiction** $A \cap A^C = \emptyset$: An element cannot simultaneously belong to $A$ and not belong to $A$. In probabilistic terms, an event cannot happen and not happen at the same time.
- **Law of the Excluded Middle** $A \cup A^C = X$: A random element $x \in X$ exclusively belongs to $A$ or $A^C$. In probabilistic terms, only one of $A$ and $A^C$ can take place.

The above led Seferis to seek the essence of whether the contradictions he often encounters can be explained by applying classical logic. In the case of fuzzy sets, the above mentioned laws do not apply.

**Fuzzy set definition:** Given $X$ is a classical reference set, every function $A:X \rightarrow [0,1]$ is called a fuzzy subset of $X$. If $x \in X$ then value $A(x)$ is called the membership value of $x$ and reverberate the membership grade of $x$ in the fuzzy set, in other words, the grade of truth of the proposition. Namely, between 0 and 1, fuzzy logic does exist.

The following notation is often used in order to define a fuzzy subset of a set $X$: $(X, \mu_A)$ where $\mu_A:X \rightarrow [0,1]$ is considered to be a fuzzy subset $A$ of $X$. For example, high temperature $A(x)$ can be depicted by applying fuzzy set theory and using the graph of the function as follows:

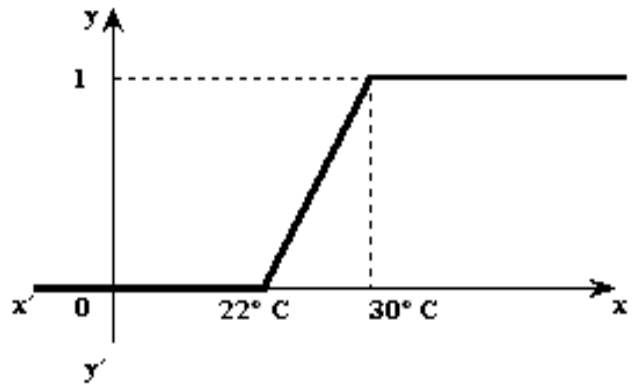

$$A(x) = \begin{cases} 0, & \alpha\nu \; x \leq 22 \\ \frac{1}{8}x - \frac{22}{8}, & \alpha\nu \; 22 \leq x \leq 30 \\ 1, & \alpha\nu \; x \geq 30 \end{cases}$$

The answer to why Seferis searches for what lies between 1 and 0 is that the law of contradiction does not apply to fuzzy logic, thus, he can explain the contradictions he encountered throughout his life.

The graph of $A(x)$ and $A^C(x)$ is shown below:

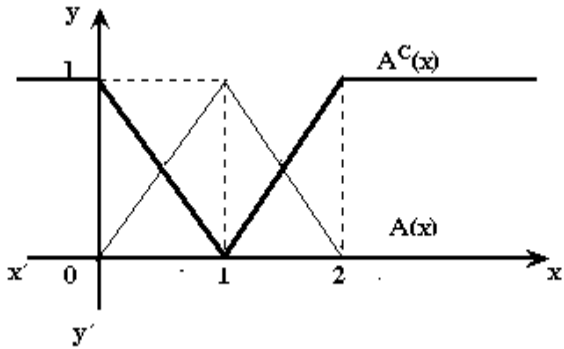

The graph of $(A \cap A^C)(x)$:

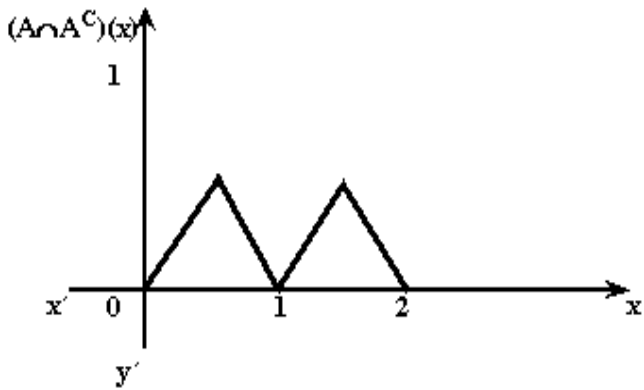

From the above graph, we gather that $A \cap A^C \neq \emptyset$. According to fuzzy set theory, $A = A^C$ can also be true at the same time. This is possible if we consider the fuzzy set: $A: X \rightarrow [0,1]$, where $A(x) = 1/2$, for all $x \in X$. In this case, we get $A = A^C$! This means we consider Seferis' "Not a God", which is the negation of 1, to be 0. This is feasible due to the above attribute of fuzzy logic negations.

At this point we should note that negation in fuzzy logic is a function $\eta: [0,1] \rightarrow [0,1]$ with the following attributes:

- $\eta(0) = 1, \eta(1) = 0$

- $\eta(\eta(x)) = x$

- It is a strictly decreasing function.

A family of negations is:

$$\eta_\lambda(x) = \frac{1-x}{1+\lambda x}, \; \lambda \in R, \; \lambda > -1$$

For $\lambda = 0$ we get the classical negation:

$\eta(x) = 1 - x$

In every case, $\eta(1) = 0$ applies and, in addition, in the case of $\eta_\lambda$ negations,

$\eta_\lambda(x) < 1 - x, \; \forall \in (0,1)$ is also possible to apply.

## 5. Comments – Conclusions.

"Helen" is one of the more mature specimens of Seferis mythical poetic era. The marvelous technique allows a poetic artifact that makes the past identify with the present of the poet. So the historical time is discharged from the subjective and social conventions and allows a unconscious parallelism and identification with mathematical concepts and sets foundations for a Fuzzy Logic. If we carefully observe our everyday life, we will notice that we ourselves experience fuzzy logic through our actions and, despite all the contradictions, we manage to achieve balance and harmony.

## 6. Acknowledgements.

The authors would like to thank Georgia Papadopoulou for her ideas and for sharing her love for Seferis and Eurippides.

* *Corresponding author..*
E-mail address: